\documentclass[11pt]{amsart}
\usepackage{amssymb}
\input xypic
\xyoption{all}
\newcommand{\BZ}{{\mathbb{Z}}}

\newcommand{\BN}{{\mathbb{N}}}

\newcommand{\BR}{{\mathbb{R}}}
\newcommand{\BC}{{\mathbb{C}}}

\newcommand{\BQ}{{\mathbb{Q}}}
\newcommand{\BP}{{\mathbb{P}}}

\newcommand{\gD}{\Delta}
\newcommand{\gd}{\delta}

\newcommand{\gc}{\gamma}

\newcommand{\gO}{\Omega}

\newcommand{\ga}{\alpha}
\newcommand{\gt}{\tau}

\newcommand{\gTh}{\Theta}
\newcommand{\gT}{\Theta}

\newcommand{\cA}{{\mathcal{A}}}

\newcommand{\cO}{{\mathcal{O}}}
\newcommand{\cP}{{\mathcal{P}}}

\newcommand{\caH}{{\mathcal{H}}}

\newcommand{\ti}[1]{\tilde{#1}}

\newcommand{\ol}[1]{\overline{#1}}

\newcommand{\Pic}{\mathrm{Pic}}

\newcommand{\sm}{\smallsetminus}
\newcommand{\Cy}[1]{\BZ/#1\BZ}
\newcommand{\hra}{\hookrightarrow}

\swapnumbers
\theoremstyle{plain}
\newtheorem{lma}{Lemma}[section]
\newtheorem{thm}[lma]{Theorem}
\newtheorem{prp}[lma]{Proposition}

\newtheorem{cor}[lma]{Corollary}

\theoremstyle{definition}
\newtheorem{prd}[lma]{Proposition-Definition}
\newtheorem{dfn}[lma]{Definition}
\newtheorem{rmr}[lma]{Remark}
\newtheorem{ntt}[lma]{Notation}
\newtheorem{exm}[lma]{Example}

\newtheorem{dsc}[lma]{}

\newcommand{\cX}{\mathcal{X}}
\newcommand\Sp{\text{Sp}}
\newcommand{\calD}{\mathcal{D}}

\newcommand{\prt}[1]{{#1}^1}

\begin{document}

\title[Intersections in the Poincar\'e bundle and the universal $\gT$]
{Some intersections in the Poincar\'e bundle and the universal theta divisor on $\ol{\cA_g}$}
\author{Samuel Grushevsky}
\address{Mathematics Department, Princeton University, Fine Hall,
Washington Road, Princeton, NJ 08544, USA}
\email{sam@math.princeton.edu}
\author{David Lehavi}
\address{Mathematics Department, University of Michigan,
2074 East Hall, 530 Church St, Ann Arbor USA}
\email{dlehavi@math.umich.edu}
\thanks{}
\begin{abstract}
We compute all the top intersection numbers of divisors on the total
space of the Poincar\'e bundle restricted to $B\times C$ (where B is
an abelian variety, and $C\subset B$ is any test curve). We use
these computations to find the class of the universal theta divisor
and $m$-theta divisor inside the universal corank 1 semiabelian
variety --- the boundary of the partial toroidal compactification of
the moduli space of abelian varieties. We give two computational
examples: we compute the boundary coefficient of the Andreotti-Mayer
divisor (computed by Mumford but in a much harder and ad hoc way),
and the analog of this for the universal $m$-theta divisor.
\end{abstract}

\maketitle

%
\section{Introduction}
%
The Poincar\'e bundle is the moduli space of degree $0$ line bundles
over an abelian variety $B$. It is a central object in research of
abelian varieties. Yet even the intersection theory on the total
space of the Poincar\'e bundle is not well understood. In this work
we consider the restriction of the Poincar\'e bundle to a test curve
$C\subset B$ times the abelian variety $B$ itself. We compute the
image of the Neron-Severi group of the total space of the Poincar\'e
bundle under this restriction, and the top intersection products of
the classes in this group.

The usefulness of this computation lies in the fact that such a test
curve $C$ is one of the two generators for the cone of curves on the
partial compactification of the moduli space of principally polarized
abelian varieties. Indeed, the universal family of ``semiabelian
varieties'', which play the same role in the moduli of abelian
varieties that the rational nodal curve plays in the moduli of
elliptic curves, is naturally the total space of the global
Poincar\'e bundle over the moduli space (see below for a more precise
description).

We are thus able to use our computations in the Poincar\'e bundle to
compute the degeneration of the universal theta divisor to the
universal semiabelian varieties, with and without level. Finally, the
knowledge of the class of the universal semiabelian theta divisor
allows us to compute some numerical invariants. We demonstrate this
with two computations: a short and straightforward computation of the
boundary coefficient of the closure of Andreotti-Mayer divisor (the
divisor consisting of principally polarized abelian varieties for
which the theta divisor is singular --- see \cite{anma}).

This divisor was computed by Mumford in \cite{mumford} by
geometrically interpreting the condition of the semiabelian theta
divisor being singular; Yoshikawa in \cite{yoshikawa} has constructed
an explicit modular form defining the Andreotti-Mayer divisor, and
thus obtained a formula for its class. The advantage of our approach
is that once the intersection-theoretic expression for the divisor on
the open part of $\cA_g$ is known, this same expression is used on
the universal semiabelian family, with no further geometric intuition
required. Using our machinery we are also able to compute the
boundary coefficient of a level-$m$ analog of the Andreotti-Mayer
divisor: the divisor on the level-$m$ cover of $\cA_g$ consisting of
principally polarized abelian varieties for which the $m$-theta
divisor is singular.

Throughout this paper we work over the complex numbers.

%
\section{Top intersections of divisors on the Poincar\'e bundle}\label{Stop}
%
\begin{ntt}[Moduli space of abelian varieties]
In this paper we consider the moduli space of principally
polarized complex abelian varieties (ppavs) of dimension $g$, which
we denote $\cA_g$. It is a stack, and we will need to be careful with
the stackiness, especially when dealing with the level cover and
branching along the boundary. However, for the Grothendieck-Riemann-Roch
computations that we do stackiness would not be a problem. We denote
by $\cX_g\to\cA_g$ the universal family of ppavs, with the fiber over
a point $[A]$ being the abelian variety $A$ itself.
\end{ntt}

\begin{ntt}
Throughout the paper, for an algebraic variety (or a Deligne-Mumford
stack) $X$ we will denote by $NS(X)$ the Neron-Severi group of
numerical equivalence classes of divisors on $X$, and by $CH^*(X)$
the Chow ring of $X$.
\end{ntt}

\begin{ntt}
For a ppav $B\in\cA_g$ we denote by $\Theta_B\subset B$ the
(symmetric) divisor of its principal polarization. For $m\in\BZ$ we
denote by $m_B:B\to B$ the multiplication by $m$ map, $m_B(z):=mz$ on
$B$. This map will turn out to be important for level cover
considerations. Finally for any point $b\in B$ we denote by $\gt_b:B\to
B$ the translation by $b$ map $\gt_b(z):=z+b$.
\end{ntt}

\subsection*{The Neron-Severi group $NS(B\times B)$ and its restriction to $NS(B\times C)$}
\begin{ntt}
Let $B\in\cA_{g-1}$ be a very general ppav, and let $C\subset B$ be a
very general curve of degree $n$ in it, i.e. such that $C\cdot
\gT_B=n$. In the following text we denote by $r$ the map $C\to B$,
and consider the image of the restriction map
\[
  r^*:NS(B\times B)\to NS(B\times C).
\]

Let $\cP$ be the Poincar\'e bundle over $B\times B$ --- it is the
universal degree $0$ line bundle over $B$, i.e. the unique line
bundle such that $\cP|_{0\times B}$ is trivial, and $\cP|_{B\times
b}$ for any $b\in B$ is the degree zero line bundle on $B$
corresponding to $b$, when we identify $B$ with the dual abelian
variety $B^{\vee}=\Pic^0(B)$ by using the principal polarization on
$B$.

Let us denote by $\ga=c_1(\cP)$ the first Chern class of the
Poincar\'e bundle on $B\times B$ (restricted to $B\times C$). Let
$\mu,\eta$ be the restrictions of $\gT_B\times B$ and $B\times\gT_B$
to $B\times C$. We also define three curves in $B\times C$:
\[
  \mu^*:=\lbrace (x,0)|x\in C\rbrace,\quad \eta^*:=\lbrace
  (0,x)|x\in  C\rbrace, \quad \gd^*:=\lbrace (x,x)|x\in C\rbrace.
\]
Finally denote by $s$ the ``shift'' automorphism
\[
  \begin{aligned}
  s:B\times B&\to B\times B\\
   (z,b)&\mapsto (z+b,b),
  \end{aligned}
\]
and note that $s$ restricts to an automorphism of $B\times C$. We
then denote by $s^*$ the action of $s$ on $NS(B\times B)$ by
pullback.
\end{ntt}
\begin{prp}\label{Pinter}
We have the following intersection numbers on $B\times C$:
\[
 \begin{array}{c|ccc}
 {}_{\rm curve}\!\!\diagdown^{\rm divisor}&\mu&\eta&\ga\\ \hline
 \mu^*&n&0&0\\
 \eta^*&0&n&0\\
 \gd^*&n&n&2n\\
 \end{array}
\]
\end{prp}
\begin{proof}
Computing the intersection of the classes $\mu$ and $\eta$ with these
curves is easy: we just forget the irrelevant factor. To compute the
intersection of the class $\ga$ of the Poincar\'e bundle with these
curves note that the Poincar\'e bundle is trivial on $0\times C$ and
$B\times 0$, so that the intersections with $\mu^*$ and $\eta^*$ are
zero, while the restriction of the Poincar\'e bundle to the diagonal,
pulled back to one of the factors, is $\cO(2\gTh_B)$ (see
\cite{mumford}, Second statement of Proposition 1.8), the degree of
which on $C$ is $2n$. Thus we are done.
\end{proof}
\begin{cor}
The group $r^*(NS_\BQ(B\times B))\subset NS_\BQ(B\times C)$ is
generated by $\ga,\mu,\eta$.
\end{cor}
\begin{proof}
By \cite{bila} $NS(B\times B)$ is $3$-dimensional, since $B$ is very
general and does not have automorphisms, but we have shown that
$\ga,\mu,$ and $\eta$ are linearly independent, since their
intersections with three curves are linearly independent.
\end{proof}
\begin{prd}
For $N\in\BZ$ let $\mu_N:=(s^*)^N(\mu)$,
then in $NS(B\times C)$ we have $\mu_N\equiv \mu+N\ga+N^2\eta$.
\end{prd}
\begin{proof}
Note that
\[
  \begin{aligned}
  s^N(\mu^*)=&\{(x,0)|x\in C\}=\mu^*,
  \quad s^N(\eta^*)=\{(Nx,x)|x\in C\},\\
  \quad s^N(\gd^*)=&\{(N+1)x,x)|x\in C\}.
 \end{aligned}
\]
We recall that if $E\in CH^e(B)$, then $[(m_B)_*E]=m^{2e}E$ (using
Poincar\'e duality and applying Theorem 6.2 from \cite{milne}). Hence
\[
  \begin{aligned}
  \mu_N\cdot\mu^*=&\mu\cdot s^N(\mu^*)=\mu\cdot\mu^*=n,\\
  \mu_N\cdot\eta^*=&\mu\cdot s^N(\eta^*)=
  \langle\gT_B\cdot\{Nx|x\in C\}\rangle_B=N^2\langle \gT_B\cdot C\rangle_B
  = N^2n\\
  \mu_N\cdot\gd^*=&\mu\cdot s^N(\gd^*)=
  \langle\gT_B\cdot\{(N+1)x|x\in C\}\rangle_B\\
  =&(N+1)^2\langle \gT\cdot C\rangle_B =(N+1)^2n.
  \end{aligned}
\]
The result now follows from the intersection matrix in Proposition
\ref{Pinter}.
\end{proof}
\begin{cor}\label{shiftact}
The action of $s^*$ on $r^*(NS(B\times B))$ is given by
\[
  s^*(\mu)=\mu+\ga+\eta;\quad s^*(\ga)=\ga+2\eta;\quad
  s^*(\eta)=\eta.
\]
\end{cor}
\begin{proof}
By the previous Proposition we have for all $N$:
\[
  \begin{aligned}
  s^*(\mu+N\ga+N^2\eta)=&s^*\mu_N
  =\mu_{N+1}=\mu+(N+1)\ga+(N+1)^2\eta\\=&(\mu+\ga+\eta)+N(\ga+2\eta)+N^2\eta.
  \end{aligned}
\]
The result follows by
equating the coefficients of the corresponding powers of $N$ on
both sides.
\end{proof}

\subsection*{Top intersections in the Chow ring of $B\times C$}
\begin{prp}
The top intersections numbers of divisors in $r^*(NS(B\times B))$ on
$B\times C$ are completely determined by the following relations:
\[
  \begin{aligned}
  (\square)& \quad  \eta^2=0
  \qquad\qquad&(\diamondsuit)\quad& \eta\mu^{g-1}&=&\ n(g-1)!\\
  (\triangle)&\quad \ga\eta=0
  \qquad\qquad&(\heartsuit)\quad& \ga^k\mu^{g-k}&=&
        \ \begin{cases}-2(g-2)!{\rm\ if\ } k=2\\
                      0{\rm\ otherwise}
        \end{cases}
  \end{aligned}
\]
\end{prp}
\begin{proof}
equation $\square$ follows from the fact that $\eta$ is a pullback of
a class on a curve.

For equation $\triangle$, note that the Poincar\'e bundle is
numerically trivial on all horizontal and vertical fibers of
$B\times B$. Since $\ga\cdot\eta$ is the $c_1$ of the restriction of
the Poincar\'e bundle to $\eta$, which is geometrically just $B$
times $n$ points, which is zero.

To prove $\diamondsuit$, we compute
\[
  \eta\mu^{g-1}=\mu^{g-1}|_\eta=\mu^{g-1}|_{B\times\lbrace n\text{ points}
  \rbrace}
  =n\gTh_B^{g-1}|_B=n(g-1)!,
\]
since restricted to $B\times pt$, the class $\mu$ is just $\gTh_B$.

To prove equation $\heartsuit$ note that the class $\mu_N$ is a
pullback of $\gTh_B$ from $B$ to $B\times C$ under the map $(z,b)\to
(z+Nb)$. Thus, as a pullback class from a $(g-1)$-dimensional
variety, its $g$'th power is zero:
\[
\begin{aligned}
  &0=\mu_N^g=((\mu+N\ga)+N^2\eta)^g\overset{\square}{=}
  (\mu+N\ga)^g+N^2g\eta(\mu+N\ga)^{g-1}\\
  &\overset{\triangle}{=}(\mu+N\ga)^g+N^2g\eta\mu^{g-1}
  \overset{\diamondsuit}{=}(\mu+N\ga)^g+N^2gn(g-1)!.
\end{aligned}
\]
This equality must hold for all $N$, which implies that in the
binomial expansion of $(\mu+N\ga)^g$ only the term quadratic in $N$
can be non-zero. Since the term quadratic in $N$ is
$N^2g(g-1)\ga^2\mu^{g-2}/2$, the equality $\heartsuit$ follows.
\end{proof}

\subsection*{Top intersections in the Chow ring of the total space
of the Poincar\'e bundle}
\begin{ntt}
Let $E$ be the trivial line bundle on $B\times B$. Denote by
$\ti{\cP}\to B\times C$ the restriction of the bundle
$\BP(E\oplus\cP)$ from $B\times B$ to $B\times C$, and set
$\xi:=c_1(\cO_{\ti{\cP}}(1))$. We define also the $0,\infty$ section
of this bundle by setting:
\[
 \cP_0:=\BP(E\oplus\{0\})|_{B\times C}\quad\cP_\infty:=\BP(\{0\}\oplus\cP)|_{B\times C}.
\]
We denote the pullback of any class $F\in CH^*(B\times C)$ to
$CH^*(\ti{\cP})$ by $\ti{F}$.
\end{ntt}
\begin{prp}\label{Pstandard}
We list some standard properties of vector bundles in our context:
\begin{enumerate}
\item The class of the vertical fiber $\BP^1$ of $\ti{\cP}\to
    B\times C$ has intersection $1$ with $\xi,$ and $0$ with
    pullbacks of all the divisors from $B\times C$.
\item\label{Iiso} The projection $\cP_0\to\BP(E)\to B\times C$ is
    an isomorphism.
\item\label{Ic1} In the Neron-Severi group $NS(\ti{\cP})$ we have
    $\cP_0-\cP_\infty\equiv\widetilde{c_1(\cP)}-\widetilde{c_1(E)}=\ti{\ga}$
\item\label{Ixi} In the Neron-Severi group $NS(\ti{\cP})$ we have
    $\cP_\infty\equiv\xi$.
\end{enumerate}
\end{prp}
\begin{prp}\label{Chowring}
The Chow ring of $\ti{\cP}$ is given by
\[
  CH^*(\ti{\cP})=CH^{*}(B\times C)[\xi]/(\xi^2+\ti{\ga}\xi).
\]
\end{prp}
\begin{proof}
Since $\ti{\cP}=\BP(E\oplus \cP)$ we have (\cite{fulton}, remark
3.2.4)
\[
  0=\xi^2+c_1(E\oplus \cP)\xi+c_2(E\oplus \cP)=\xi^2+\ti{\ga}\xi.
\]
\end{proof}
\begin{rmr}
Note that this matches with the obvious fact that the 0 and infinity
sections do not intersect: indeed, we have
$\cP_0=\cP_\infty+\ti{\ga}=\xi+\ti{\ga}$, and so it checks out that
\[
 \cP_0\cP_\infty=(\xi+\ti{\ga})\xi=0.
\]
\end{rmr}
\begin{prp}
The top intersection products of divisors on $\ti{\cP}$ are
completely determined by the pullback of the relations
$\square,\triangle,\diamondsuit,\heartsuit$ together with:
\[
  \begin{aligned}
  (\blacksquare)\quad \xi^2=-\ti{\ga}\xi,\qquad
  (\blacktriangle)&\ \begin{matrix}  \text{ the map }
  CH^*(B\times C)\to CH^*(\xi)\text{ arising from}\\
 \text{statements \ref{Iiso},\ref{Ixi} in Prop.
 \ref{Pstandard} is an isomorphism},\end{matrix}
  \end{aligned}
\]
and the fact that top intersections of pullback classes are $0$ (for
dimension reasons).
\end{prp}
\begin{rmr}[Motivation for the notations]
We only consider intersections of two divisors or top intersections.
The identities for top intersection numbers are denoted by card
suits, while the equations for the squares of the generators are
denoted by squares. The triangles thus denote relations for products
of two different classes. Moreover, the  pullback relations are
white, while the relations involving the fiber generator $\xi$ are
black.

This means that in a computation of a top intersection number on
$\ti{\cP}$ we would typically do the following: first apply black
relations (and the fact that the top intersection of pullback classes
is zero) to reduce the computation to a computation on $B\times C$,
then apply $\square$ and $\triangle$ to get rid of redundant
intersections, and finally the card suit relations to get actual
numbers.
\end{rmr}
\begin{rmr}[Computational Trick]
We conclude this section with a small computational trick, coming
from the relations we already deduced:
\[
\begin{aligned}
  &(\xi+a\ti{\mu}+b\ti{\ga})^{g+1}=\sum\limits_{k=0}^{g+1}
  \binom{g+1}{k}\xi^k(a\ti{\mu}+b\ti{\ga})^{g+1-k}\\
  \mathop{=}\limits^{\blacksquare}
  &\sum\limits_{k=1}^{g+1}\binom{g+1}{k}\xi(-\ti{\ga})^{k-1}(a\ti{\mu}+b\ti{\ga})^{g+1-k}
  \mathop{=}\limits^{\blacktriangle}\sum\limits_{k=1}^{g+1}\binom{g+1}{k}(-\ga)^{k-1}(a\mu+b\ga)^{g+1-k}\\
  \mathop{=}\limits^{\heartsuit}&\sum\limits_{k=1}^3\binom{g+1}{k}(-\ga)^{k-1}
  \binom{g+1-k}{3-k}(b\ga)^{3-k}(a\mu)^{g-2}\\
  =&\binom{g+1}{3}(a\mu)^{g-2}\ga^2(1-3b+3b^2)\\
\mathop{=}\limits^{\heartsuit}&-n\binom{g+1}{3}a^{g-2}2(g-2)!(b^3-(b-1)^3)=-\frac{n(g+1)!}{3}a^{g-2}(b^3-(b-1)^3).
  \end{aligned}
\]
We will denote this relation in the sequel by ``$T$'' (for
``trick'').
\end{rmr}
%
\section{Ppavs and rank one degenerations}\label{Snolev}
%
\subsection*{Rank 1 degenerations of abelian varieties, and moduli}
\begin{dfn}[Semiabelian varieties]\label{semiabelian}
A (non-normal compactification of a rank-one-degenerated principally
polarized complex) semiabelian variety is constructed as follows. Let
$B\in\cA_{g-1}$, $b\in B$, and let $S\to B$ be the line bundle
corresponding to $b$ under the identification $B\cong\Pic^0 B$.
Complete $S$ fiberwise to a $\BP^1$-bundle $\ti{S}\to B$, and then
identify the $0$ and $\infty$ sections of $\ti{S}$ (each a copy of
$B$) with a shift by $b$, i.e. define $\ol{S}:=\ti{S}/(z,0)
\sim (\gt_bz,\infty)$ (recall that we defined $\gt_b$ to be the map
$z\mapsto z+b$ on $B$).
The divisor of the principal polarization
$\gTh\subset\ol{S}$ is a section $\left(B\sm(\gTh_B\cap
\gt_b\gTh_B)\right)\to \ti{S}$, together with the entire fibers over
$\gTh_B\cap \gt_b\gTh_B$, such that it intersects the $0$-section of
$\ti{S}$ in $\gTh_B$, and correspondingly the $\infty$-section of
$\ti{S}$ in $\gt_b\gTh_B$, thus giving a well-defined divisor on
$\ol{S}$. Notice that the compactification $\ol{S}$ no longer admits
a projection map to $B$. We also note that if we started with the
point $-b$ instead of $b$, the resulting object is going to be
isomorphic --- indeed, any semiabelian variety has an involution
interchanging the 0 and $\infty$ sections, and this involution sends
a semiabelian variety with bundle $b$ on the base (the 0 section) to
the one with bundle $-b$. This means that the family of semiabelian
varieties with base $B$ is parameterized by $B/\pm1$, and thus this
family is singular at points of order two on $B$.
\end{dfn}
\begin{dfn}[Partial compactification]\label{Dprt}
We denote $\prt{\cA_g}$ the partial compactification of $\cA_g$
obtained by ``adding'' the locus of rank 1 semiabelian varieties: See
\cite{Amoduli} (setup 1.2.8), \cite{olsson} for the construction of
the second Voronoi toroidal compactification $\overline{\cA_g}$ of
$\cA_g$, over which there exists a universal family
$\overline{\cX_g}$. The partial compactification  $\prt{\cA_g}$ is a
subscheme of $\overline{\cA_g}$.

The locus of rank 1 semiabelian varieties in $\prt{\cA_g}$ forms the
boundary divisor denoted $\gD\subset\prt{\cA_g}$. We denote by
$\prt{\cX_g}$ the partial compactification of the universal family of
abelian varieties, which is an extension of the family
$\cX_g\to\cA_g$ such that the fibers over the boundary are
semiabelian varieties. We denote by $\calD\subset\prt{\cX_g}$ the
universal divisor of the principal polarization, i.e. the universal
theta divisor.
\end{dfn}
\begin{dfn}[The universal rank 1 degeneration --- see \cite{mumford} 1.8]
The universal family of semiabelian varieties in $\partial
\prt{\cX_g}$ lying over $[B]\in\cA_{g-1}$, i.e. the fiber over $[B]$
of the map $\prt{\cX_g}\to\prt{\cA_g}\to\cA_{g-1}$, is defined in the
following way: consider the $\BP^1$-bundle $\ti{\cP}=\BP(E\oplus
\cP)$ over $B\times B$. Denote by $\cP_0$ and $\cP_\infty$ the zero
and infinity sections of $\ti{\cP}$ over $B\times B$, respectively.
Then the universal semiabelian variety over $B\times B$ (without
level structure) is obtained from $\ti{\cP}$ by identifying $\cP_0$
and $\cP_\infty$ with a shift by $b$:
\[
  \ol{\cP}:=\ti{\cP}/(z,b,0)\sim(s(z,b),\infty),
\]
and further taking the quotient under $b\mapsto -b$.
\end{dfn}
From now on we will restrict this construction to our test curves,
and use the same letters to denote the restrictions of the glued
bundle and its normalization from $B\times B$ to $B\times C$. Since
we chose the curve $C\subset B$, we can make sure that $C\cap
(-C)=\emptyset$, and thus in working with $B\times C$ we can ignore
the sign involution of $b$.

\begin{rmr}
It is known since Tai's work in \cite{tai} that one can construct
some $\cA_g'\supset\prt{\cA_g}$ such that $\cA_g'\sm\prt{\cA_g}$ has
codimension $2$ in $\cA_g'$, and $\cA_g'$ has only canonical
singularities (the rigorous proof of this is given in
Shepherd-Barron's \cite{shepherdbarron} and its addendum, to appear).
\end{rmr}

\subsection*{The universal theta divisor on rank $1$ degenerations}

\begin{prp}\label{Ptoptheta}
The top intersection number of the universal theta divisor over
$B\times C$ restricted to $\cP_0$ is $0$.
\end{prp}
\begin{proof}
By definition on one semiabelian variety $\ti{S}\to B$ the
restriction of the theta divisor to the 0-section
$B=S_0\subset\ti{S}$ is just $\Theta_B$, and thus
$\calD|_{\cP_0}=p^*\Theta_B$, where $p:B\times C\to B$ is the
projection. Hence the top intersection is a pullback of a top+$1$
intersection product on $B$.
\end{proof}
\begin{thm}\label{Tthetaclass}
The class of the universal theta divisor on the normalization
$\ti{\cP}$ of the universal semiabelian variety, restricted to
$B\times C$, is equal to
\[
  \calD=\xi+\ti{\mu}+\frac{1}{2}\ti{\ga}+\frac{1}{4}\ti{\eta}
\]
in $r^*(NS(\ti{\cP}))$.
\end{thm}
\begin{rmr}
Below we give an elementary algebraic derivation of this formula from
the geometric description of the universal semiabelian family.
Alternatively one could work out explicitly the degenerations of the
theta function along the boundary --- this approach is followed in
\cite{huwe} for the case of $g=2$, and one easily computes the class
of the universal semiabelian theta function given there to be what we
claim. Our approach has the advantage of giving a way to
straightforwardly deal with the computation on the level cover as
well.
\end{rmr}
\begin{proof}
Denote the coefficients
\[
  \calD=c_\xi\xi+c_\mu\ti\mu+c_\ga\ti\ga+c_\eta\ti\eta.
\]
We first compute $c_\xi$; this coefficient can be computed after
intersecting with $\eta$, i.e. on the fiber over one moduli point in
$\partial \cA_g$. We compute it by considering two families in
$\prt{\cA_g}$ with the same flat limit, which is the trivial bundle
over $B$ with 0 and $\infty$ sections identified by identity. Our
first family is the family of semiabelian varieties parameterized by
$(B,b)$ where $b\to 0$. Our second family is the family $B\times E$
where $E$ is a moving elliptic curve that degenerates to the rational
nodal curve. Since the theta divisor on the second family is
numerically given by $\gT_B\times E+B\times pt$, the limit theta
divisor is given by $\gT_B\times(\BP^1/0\sim\infty)+B\times pt$.
Since the intersection of this limit theta divisor with a general
fiber of the bundle is $1$ (coming from the $B\times pt$ component),
the same holds in the first family.

We note that with a little more work we could have recovered the
coefficient of $\mu$, but we will get it almost for free below.

Next we compute the coefficients $c_\mu,c_\ga$: note that by
($\blacksquare$)
\[
 \calD|_{\cP_\infty}=(\xi+c_\mu\ti{\mu}+c_\ga\ti{\ga}+c_\eta\ti{\eta})|_\xi=
 (c_\mu\ti{\mu}+(c_\ga-1)\ti{\ga}+c_\eta\ti{\eta})|_{\cP_\infty}
\]
and
\[
 \calD|_{\cP_0}=(\xi+c_\mu\ti{\mu}+c_\ga\ti{\ga}+c_\eta\ti{\eta})|_{\xi+\cP}=
 (c_\mu\ti{\mu}+c_\ga\ti{\ga}+c_\eta\ti{\eta})|_{\cP_0}.
\]
Since the theta divisor is a stable limit of smooth connected
varieties, it is connected on codimension $1$. This means that if we
identify both $\cP_0$ and $\cP_\infty$ with $B\times C$, then
$\calD|_{\cP_0}$ is the pullback of  $\calD|_{\cP_\infty}$ under the
shift $s$, i.e.
\[
  \begin{aligned}
    c_\mu\mu+c_\ga\ga+c_\eta\eta=&s^*(c_\mu\mu+(c_\ga-1)\ga+c_\eta\eta)\\
  =&c_\mu(\mu+\ga+\eta)+(c_\ga-1)(\ga+2\eta)+c_\eta\eta\\
  =&c_\mu\mu+(c_\mu+c_\ga-1)\ga+(c_\mu+2c_\ga+c_\eta-2)\eta.
  \end{aligned}
\]
Comparing the $\ga$ coefficient we have $c_\mu=1$, and comparing
$\eta$ coefficients we get
\[
  c_\ga=(2-c_\mu)/2=1/2.
\]
To compute $c_\eta$, we observe that by Proposition \ref{Ptoptheta}
we have
\[
\begin{aligned}
 0=&(\calD|_{\cP_\infty})^g
  \overset{\blacksquare,\blacktriangle}{=}
  (\mu-\frac{1}{2}\ga+c_\eta\eta)^g
  \overset{\square,\triangle}{=}gc_\eta\eta\mu^{g-1}
  +(\mu-\frac{1}{2}\ga)^g\\
  \overset{\heartsuit,\diamondsuit}{=}&gc_\eta(g-1)!+
\binom{g}{2}(-2(g-2)!)\frac{1}{4}=g!\left(c_\eta-\frac{1}{4}\right)
\end{aligned}
\]
\end{proof}
\begin{dsc}
In \cite{mumford} Mumford computed the class in
$NS(\prt{\cA_g})$ of the closure in $\prt{\cA_g}$ of the
Andreotti-Mayer divisor. Since
$NS(\prt{\cA_g})$ is spanned by the Hodge class and the boundary class,
he had to compute these two coefficients.

Mumford's computation of the Hodge coefficient of the class of the
Andreotti-Mayer divisor (in \cite{mumford} Proposition 3.1) uses the
ramification formula (see \cite{fulton} Example 9.3.12):
\[
  \BR(f)=(c(f^*(T_C))c(T_\calD)^{-1})_g,
\]
where $f:\calD\to C$ is the universal theta divisor over a smooth
curve in $\cA_g$. In practice, Mumford breaks the morphism $f$ to a
composition $\calD\overset{i}{\hra} \cX\to C$, where $\cX$ is the
universal abelian variety over $C$, and gets (we omit obvious
pullback notations):
\[
  \begin{aligned}
  c(T_\calD)=&c(T_{\calD/C})c(T_C)\quad\Rightarrow\quad
  \BR(f)=(c(T_{\calD/C})^{-1})_g,\\
  c(T_{\calD/C})=&c(T_{\calD/\cX})c(T_{\cX/C})=(1-\calD)c(T_{\cX/C})\quad\Rightarrow\\
  i_*\BR(f)=&(c(\gO_{\cX/C})(1-\calD)^{-1})_g\cdot D
   =\sum_{i=0}^gc_i(\gO_{\cX/C})\calD^{g+1-i}.
  \end{aligned}
\]
Mumford proceeds to compute the class --- see the second part of
Proposition 3.1 in \cite{mumford}.

To compute the boundary coefficient, Mumford uses two hard
ad-hoc arguments. However, armed with our understanding of the theta divisor
we can give an alternative straightforward computation:
\end{dsc}
\begin{prp}\label{Pmumbound}
Let $C\subset B$ be a test curve as above, let $f:\ol{\calD}\to C$ be
the universal semiabelian theta divisor over the test curve $C$ (so
that its lifting to $\ti{\cP}$ from $\ol{\cP}$ is $\calD$); then
$\BR(f)=\frac{n(g+1)!}{6}$.
\end{prp}
\begin{proof}
Since the divisor $\calD$ is singular we cannot use the ``smooth''
ramification formula. Instead we use the ``singular'' version:
\[
  \BR(f)=(c(f^*(T_C))s(\gD_{\ol{\calD}},\ol{\calD}\times\ol{\calD}))_g.
\]
See \cite{fulton} the end of example 9.3.12 for details (compare
Johnson/Fulton-Laksov singular double point formula in 9.3.6 to the
double point formula in Theorem 9.3).

We compute this intersection product by pulling back under
$\pi:\calD\to\ol{\calD}$. By the definition of the Segre class we
have:
\[
  \begin{aligned}
  \pi^*\BR(f)=&(c(\pi^*f^*(C))
  s((\pi\times\pi)^*\gD_{\ol{\calD}},\calD\times\calD))_{2g-1}\\
  =&(c(\pi^*f^*(C))c(T_{\calD/\ol{\calD}})c(\calD)^{-1})_g.
  \end{aligned}
\]
By definition we have
$c(T_{\calD/\ol{\calD}})=(\cP_0+\cP_\infty)|_\calD=(2\xi+\ti\ga)|_\calD$.
However, by example \cite{fulton} 3.2.11, the total Chern class of
the relative tangent bundle is
\[
 \begin{aligned}
  c_t(T_{\ti{\cP}/C)} =&(1+\xi t)^2c_{\frac{t}{1+\xi t}}(E \oplus \cP)\\
=&(1+\xi t)^2c_{t(1-\xi t)}(\cP)=(1+2\xi t+\xi^2t^2)(1+\ti\ga t(1-\xi  t))\\
=&1+(2\xi+\ti\ga)t+(\xi^2-2\ti\ga\xi+2\xi\ti\ga)t^2
\overset{\blacktriangle}{=}1+(2\xi+\ti\ga)t.
\end{aligned}
\]
Hence, factoring the map $\calD\to C$ through
$\calD\overset{i}{\hra}\ti{\cP}\to C$, similarly to the case of
smooth abelian varieties, we have
\[
  i_*\pi^*\BR(f)=\left[\sum c_i(1+2\xi+\ti{\ga}-(2\xi+\ti{\ga}))D^{g-i}\right]_{\deg g}\cdot\calD
\]
Hence we have:
\[
\begin{aligned}
 i_*\pi^*\BR(f)=&\calD^{g+1}=\left(\xi+\ti{\mu}+\frac{1}{2}\ti{\ga}+\frac{1}{4}\ti{\eta}\right)^{g+1}\\
=&(g+1)\frac{\ti{\eta}}{4}\left(\xi+\ti{\mu}+\frac{\ti{\ga}}{2}\right)^g
 +\left(\xi+\ti{\mu}+\frac{\ti{\ga}}{2}\right)^{g+1}\\
 \overset{T}{=}&(g+1)\frac{\ti{\eta}}{4} g\xi{\ti{\mu}}^{g-1}
 +\frac{n(g+1)!}{3}\left(\left(\frac{1}{2}\right)^3-\left(\frac{-1}{2}\right)^3\right)\\
 =&n(g+1)!\left(\frac{1}{4}-\frac{1}{12}\right)=\frac{n(g+1)!}{6}.
\end{aligned}
\]
\end{proof}
%
\section{Geometry of the level cover $\cA_g(m)$ of $\cA_g$}
%
\subsection*{Background: level structure on rank 1 degenerations of abelian varieties}
\begin{dfn}[$m$-torsion points and level moduli space]
Given a ppav $A\in\cA_g$ and $m\in\BN$, we denote by $A[m]$ the set
of points of order $m$ in $A$; note that since $m$ is prime to the
characteristic of the base field, which we assumed to be $0$, we have
$A[m]\cong(\Cy{m})^{2g}$.

We denote by $\cA_g(m)$ the moduli space of ppavs together with a
choice of a symplectic basis for the set of $m$-torsion points, up to
isomorphisms. Her e symplectic means with respect to the Weil pairing
on $A[m]$ --- see \cite{milne} I.12 for a complete definition. The
precise definition of the pairing would not matter for us, and thus
we do not give it.

The forgetful map $\cA_g(m)\to\cA_g$ is a finite cover with fiber
over $[A]\in\cA_g$ being the set of all symplectic bases for
$A[m]\cong (\Cy{m})^{2g}$, and with the deck group  $\Sp(2g,\Cy{m})$.

The fine moduli stack $\cA_g(m)$ carries a universal family
$\cX_g(m)\to\cA_g(m)$, where the forgetful map lifts to the
multiplication by $m$ map $\left[A,\substack{\text{basis}\\\text{for
}A[m]}\right]\overset{m_A}{\to} A/A[m]$ over every moduli point
$[A]\in\cA_g$.
\end{dfn}
\begin{dfn}[Global theta divisors]
We define the universal order $m$ theta divisor on $\cX_g(m)$ to be
the pullback of the universal theta divisor $\calD$ under the map
\[
  \begin{aligned}
  q:\cX_g(m)&\longrightarrow\cX_g\\
[A,\{a_1,\ldots a_g,a_1',\ldots a_g'\}]&
  \mathop{\longmapsto}\limits^{id}[A/\mathrm{span}(a_1,\ldots, a_g)],
 \end{aligned}
\]
where $a_1,\ldots a_g,a_1',\ldots a_g'$ denotes a symplectic basis of
$A[m]$.
\end{dfn}
\begin{dsc}
Below we explain in what sense we have a level structure on a
semiabelian variety. We do this by analyzing  the limit of the map
$q:\cX_g(m)\to\cX_g$ defined above. We now fix (up until Proposition
\ref{Plimitq}) a moduli point $[B]\in \cA_{g-1}\subset
\cA_g^{\mathrm{Sat}}$, and consider the limit of the map $q$ over the
preimage of $[B]\in\partial\prt{\cA_g}=\cX_{g-1}/\pm 1$.
\end{dsc}
\begin{dfn}[Level $m$ rank-$1$ semiabelian varieties]\label{Dlevel}
The definition below is the rank $1$ case of the more general
definition given in \cite{Amoduli} set-up 1.2.8.

Let $b$ be a non zero point in $B$, and let $S$ be the degree $0$
line bundle over $B$ associated to $b$. Let then
$\ti{S}^{(i)}:=\BP(E\oplus S)$, for $i=0,\ldots m-1$, be $m$ copies
of the projectivization of this bundle, let $S^{(i)}$ be
$\ti{S}^{(i)}$ minus the $0$ and $\infty$ sections, and let
$\ti{S}(m):=\sqcup\lim\limits_{i=0}^{m-1} \ti{S}_i$ be the disjoint
union of $m$ copies. The (non-normal compactified rank 1) level $m$
semiabelian variety is
\[
 \ol{S}(m)=\ti{S}(m)/(z,0)^{(i+1)}\sim(z,\infty)^{(i)}\ {\rm and} \ (z,
 0)^0\sim(z+b,\infty)^{(m-1)}.
\]
See example 4a in page 40 of \cite{huleknotes} for a ``graphical''
description.

We denote by $\prt{\cA_g(m)}$ the partial compactification of
$\cA_g(m)$ obtained by ``adding'' the locus $\gD(m)$ of level $m$
rank $1$ semiabelian varieties, in the same sense of Definition
\ref{Dprt}: By \cite{Amoduli} setup 1.2.8, $\prt{\cA_g(m)}$ is the
underlying scheme of a substack of projective fine moduli stack: the
second Voronoi compactification of $\cA_g(m)$. The fibers of the
universal family $\cX_g(m)^1$ over $\prt{\cA_g(m)}$ are as defined
above.
\end{dfn}
\begin{dsc}
We note that no claim was made (yet) on how the level semiabelian
varieties form a universal family, i.e. on how the semiabelian
varieties with $B$ and $b$ varying fit together in a family; this
follows below.
\end{dsc}
\begin{dsc}[Local monodromy]\label{Fmono}
Let $\pi:U\to M$ be a fine moduli stack, and let $\gc$ be an
automorphism of $M$ with a fixed point $p$. Then $\gc$ induces an
automorphism of the fiber over $p$.
\end{dsc}
\begin{lma}\label{Lmono}
Let $[B,b]\in\partial\cA_g(m)$ be a boundary point, then the
automorphism induced on $\sqcup_{i=0}^{m-1} \ti{S}_i/B$ from the
local monodromy (i.e. by an element of $\Sp(2g,\Cy{m})$ fixing the
point $[B,b]$) in the sense of \ref{Fmono} permutes cyclically the
$m$ components of a semiabelian variety.
\end{lma}
\begin{proof}
Recall that the local monodromy about any point in a fine moduli
space is finite. Since the monodromy is discrete we may assume that
$b,B$ are general. We argue by contradiction. Consider a morphism
which fixes the components above, then on each component the
automorphism lifts to an automorphism of the fiber bundle
corresponding to $b$ over $B$. Assume that an induced automorphism
fixes the components, then it also sends the $0$ and  $\infty$
sections of each component (which are the intersections of two
components, and thus fixed) to themselves. For $B$ general the group
of finite order automorphisms of $B$ is generated by the Kummer
involution and translates by torsion points. However the torsion
points of $B$ are limits of torsion points coming from smooth abelian
varieties, and the Kummer involution is the limit of Kummer
involutions on smooth abelian varieties.

It remains to prove that if just one component is fixed then its
$0,\infty$ sections are fixed. We will assume the contrary and derive
an intersection theoretic contradiction: assume that the $0,\infty$
section are interchanged, then there exists an automorphism of the
Chow group of the universal $\BP^1$ bundle over $B$ (over some test
curve $C\subset B$) which permutes the $0,\infty$ sections. However,
the square of the infinity section is the first Chern class of the
bundle, while the square of the $0$ section is minus this first Chern
class. As one of these is effective and one is anti-effective (as
divisors on the $0$ and $\infty$ sections respectively), this is
impossible.
\end{proof}
\begin{prd}\label{Plimitq}
There is a commutative diagram of finite quotients and inclusions:
\[
  \xymatrix{
  \cX_g(m)\ar@{^{(}->}[d]\ar^q[r]&\cX_g\ar@{^{(}->}[d]\\
  \cX_g(m)^1\ar^{q^1}[r]&\cX_g^1.
 }
\]
Fixing a moduli point in $\cA_g^1(m)$ and working over its pullback
to $\cX_g^1(m)$, the map $q^1$ is defined on
semiabelian varieties as a quotient by
\begin{itemize}
\item limits of translation by the universal $m$ torsion points
    $\{a_1,\ldots a_g\}$, which are the $m$ torsion point
    subgroup of $B$, and
\item the component shift automorphism described in Lemma
    \ref{Lmono} above. We denote this component shift operator by
    $S$, for ``shift''. Its relationship with the automorphism
    $s$ of $B\times B$ will be made clear in what follows.
\end{itemize}
and as $q$ on abelian varieties.
\end{prd}
\begin{proof}
By Lemma \ref{Lmono} our claim follows locally about a generic
boundary point in $\cX^1(m)$. However our definition of the
automorphism group above is a global algebraic definition. Hence
$q^1$ is a quotient by a finite group of algebraic automorphism,
which factors through the fibers of the map
$\cX_g^1(m)\to\cA_g^1(m)$, as described locally in Lemma. Finally
recall that the quotient of a flat family by a finite group which
acts fiberwise is a flat family (we need the flatness in order to
have finite group quotients of projective varieties - flat families
live in ambient projective spaces); hence the quotient by $q^1$ is
$\cX^1$.
\end{proof}
\begin{rmr}
One can alternatively describe the local picture of the action of
$\Sp(2g,\Cy{m})$ on the universal family $\cX_g^1(m)$ near a boundary
point of $\cA_g^1(m)$ analytically in coordinates, by studying the
action of $\Sp(2g,\BZ)$ on the universal cover $\caH_g\times\BC^g$
(where $\caH_g$ stands for the Siegel upper half-space) of the
universal family, similar to the discussion in \cite{hulekbook}. In
this case if the degeneration corresponds to $ma_1$ (which is a
period vector of the abelian variety) going to infinity, then the
limits of points $a_2,\ldots,a_g,a_2',\ldots,a_g'$ are $m$-torsion
points on $B$ --- and thus induce automorphisms of it by translation
--- while the limit of the translation by $a_1$ is the shift $S$.
\end{rmr}
\begin{cor}
The pullback under $q^1$ of the universal theta divisor in $\cX_g^1$
to $\cX_g^1(m)$ is a flat family.
\end{cor}
\begin{proof}
Both the pullback of the universal theta divisor under the finite
group quotient $q^1$, and the universal $m$-theta divisor in
$\cX_g^1(m)$ are flat families. Since they are equal on $\cX_g(m)$,
they are equal on $\cX_g^1(m)$.
\end{proof}
\begin{cor}\label{limitq}
Let $[B,b]$ be a moduli point in $\prt\cA_g$, let $F$ be a maximal
isotropic group on $B$, considered as a ``limit'' of maximal
isotropic groups over the open part in the sense of Proposition
\ref{Plimitq},  and let $\cP'$ be the pullback of the Poincar\'e
bundle on $B/F\times B/F$ to $B\times B$ under the quotient map
$B\times B\to (B/F)\times(B/F)$. Then the restriction of the map
$q^1$ is the ``clutching'' (gluing the $\infty$-section of each
$\ti{\cP'}^{(i)}$ in some way to the $0$-section of
$\ti{\cP'}^{(i+1)}$, and the $\infty$ of $\ti{\cP'}^{(m-1)}$ to the
$0$ of $\ti{\cP'}^{(0)}$ in some way) of the morphism given in the
diagram:
\[
  \begin{aligned}
  (\sqcup_{i=0}^{m-1} \ti{\cP'}^{(i)}/B\times B)/
  \substack{\text{identification of }\\\text{Def. \ref{Dlevel} on
  fibers}}&\to \ti{\cP}/(B/F)\times (B/F)\\
  (x_i,b_1,b_2)&\mapsto (x,b_1+F,b_2+F).
  \end{aligned}.
\]
Specifically, the universal level abelian variety over $[B]$ is the
``clutching'' of $\sqcup_{i=0}^{m-1} \ti{\cP'}_i$ along the
$0,\infty$ fibers.
\end{cor}
\begin{proof}
Let $b\in B$ be a point in the ``moduli'', or right, copy of $B$. The
quotient of the fiber $\ti{S}(m)$ (corresponding to $b$) over
$B\times\{b\}\subset B\times B$ by the group describe in
\ref{Plimitq} is indeed the semiabelian variety corresponding to
$[B/F,b+F]$. Hence our statement is true when restricted to
$B\times\{b\}$.

This means that we have a map from any fiber of each $\ti{\cP'}_i$ to
the corresponding fiber of $\ti{\cP}$ defined naturally, and since a
linear map of line bundles is a bundle morphism, the universal line
bundle over $B\times B$ is the pullback of the Poincar\'e bundle from
$(B/F)\times(B/F)$.
\end{proof}
\begin{dsc}
The above corollary describes the structure of the map $q^1$ on the
universal semiabelian family over $B$, and by working carefully with
this description one can determine the appropriate gluings of the 0
and $\infty$ sections of consecutive bundles, and to compute the
classes of the bundles $\ti{\cP'}_i$.

However, for our computations it is much easier to change
coordinates: indeed, we can translate each of $B\times B$ by some
point and pull back the bundles, before applying the map $q^1$. Note
that shifting the origin on $B\times B$ (i.e. apply a translation to
it) changes the choice of the basis for $NS_\BQ(B\times B)$, as the
Poincar\'e bundle is not mapped to itself by a translation. The
choice of translations that we make is as follows: it is the one
induced by the component-shift action of $S$ (which satisfies
$q^1\circ S=q^1$, this is true over smooth abelian varieties, and the
limit is flat).
\end{dsc}
\begin{dfn}
To define our new coordinates on each of $\ti{\cP'}_0^{(i)}$, we do
not apply a translation on the bundle $\ti{\cP'}^{(0)}$, while for
any $i>0$ we take the point $S^i(0)$ to be the origin on
$\ti{\cP'}_0^{(i)}$. This corresponds to shifting coordinates on the
base of each of the bundles. We choose a basis for each
$NS(\ti{\cP'}^{(i)})$ in these coordinates as before, using
$\ga'^{(i)}:=c_1(\cP'^{(i)})$ instead of $\ga$. Since $S(\cP'^{(i)})=
\cP'^{(i+1)}$, in our new coordinates $\ga'^{(i)}=\ga'^{(i+1)}$ (as
classes in $NS(B\times B)$), while if we did not shift the origins of
each of $\ti{\cP'}_0^{(i)}$, this would not be the case.
\end{dfn}

\subsection*{The universal order $m$ theta divisor on rank 1 degenerations}
\begin{ntt}
Given a class $x\in NS(B\times B=NS(\ti{\cP'}^{(i)}_0)$ we denote by
$\ti{x}_i$ the class of its pullback to $NS(\ti{\cP'}^{(i)})$ (or
simply $\ti{x}$, if the choice of $i$ is clear). We also denote by
$\xi_i\in NS(\ti{\cP'}^{(i)})$ the first Chern class of the
tautological line bundle on the $\BP^1$-bundle $\ti{\cP'}^{(i)}$ over
$B\times B$.  As in section \ref{Stop}, we work with the image of
$r:NS(B\times B)\to NS(B\times C)$.
\end{ntt}
\begin{dsc}\label{levelbundle}
Proposition \ref{Chowring} applies to each of the bundles
$\ti{\cP'}^{(i)}$ to yield $CH^*(\ti{\cP'}^{(i)})=CH^*(B\times
B)[\xi_i]/ (\xi_i^2+\widetilde{\ga'}\xi_i=0)$. The discussion
following Proposition \ref{Chowring} also applies to each of the
$\ti{\cP'}^{(i)}$, so that we get ${\cP'}_\infty^{(i)}=\xi_i$, and
${\cP'}_0^{(i)}=\xi_i+\ti{\ga'}$.
\end{dsc}
\begin{thm}\label{Ttheta}
The class of the restriction of the universal theta divisor to
$\ti{\cP'}^{(i)}$ in $NS(\ti{\cP'}^{(i)})$ is numerically equivalent
to:
\[
  \calD_i\equiv\xi_i+m\ti{\mu}_i+\frac{1}{2}\ti{\ga'}_i+\frac{m}{4}\ti{\eta}_i
\]
Moreover, the class $\ga'$ is equal to $m\ga$, and thus we can write
\[
  \calD_i\equiv\xi_i+m\ti{\mu}_i+\frac{m}{2}\ti{\ga}_i+\frac{m}{4}\ti{\eta}_i
\]
\end{thm}
\begin{proof}
The proof is very similar to the proof of Theorem \ref{Tthetaclass},
and is very easy in our coordinates. Indeed, since $q^1\circ S=q^1$,
we have $S^*(D_{i+1})=D_i$. Since the coordinates are induced by the
action of $S$, this means that in the expression
\[
  D_i=c_\xi^{(i)} \xi_i+c_\mu^{(i)}\ti{\mu}_i+c_\ga^{(i)}\ti{\ga'}_i
  +c_\eta^{(i)}\ti\eta_i
\]
all coefficients are independent of $i$ (since for the basis we have
$S^*(\mu_{i+1})=\mu_i$, etc.).

We will thus compute this coefficients on the bundle
$\ti{\cP'}^{(0)}$, which is mapped to $\ti{\cP}$ simply by the map
$q^1$, with no shift involved. From Proposition \ref{Plimitq} it
follows that $c_\xi=1$. By Corollary \ref{limitq} we know that the
base $\ti{\cP'}^{(i)}_0=B\times B$ is mapped to
$\ti{\cP'}_0=B/F\times B/F$ simply by taking the quotient, i.e. by
the map $q\times q$ for the $g-1$-dimensional abelian varieties.

To compute the pullbacks of the classes $\mu,\ga,\eta$ under the map
$q$, note that the map $q\times q$ induces the map $q$ on the
vertical and horizontal copies of $B$ in $B\times B$, as well on the
diagonal, and the test curves sit in these abelian varieties, they
are all mapped to the image by the map $q$. Since the pullback of the
theta divisor under $q$ is a section of $m\Theta$ (this is a
classical theta function of order $m$), this means that the divisor
classes $\mu,\ga,$ and $\eta$ pull back to $m\mu,m\ga$, and $m\eta$
respectively.

Moreover, the bundle $\ti{\cP'}^{(0)}$ is the pullback of $\ti{\cP}$
under the map $q$, and thus we must have $\ga'=q^*\ga$, which by the
above implies $\ga'=m\ga$.
\end{proof}
\begin{exm}
As an example of using the machinery we will compute the boundary
coefficient in $NS(\prt{\cA_g})$ of the branch divisor of
$f:m\gT\hookrightarrow\prt{\cX_g}(m)\to\prt{\cA_g}(m)\to\prt{\cA_g}$.
Since $\prt{\cA_g(m)}\to\prt{\cA_g}$ is ramified to order $m$ over
the boundary, and since $m\gT$ is the pullback under the map $q$
above of the divisor $\gT\subset\cX_g^1$, by restricting this picture
to the test curve and computing in the universal family over it, we
see that this coefficient is simply $m^{g+1}$ times the coefficient
computed in \ref{Pmumbound}. To demonstrate our techniques we will
compute this using the method of \ref{Pmumbound}. This means we will
take a test curve $C$ in the boundary of $\prt{\cA_g}$ contracted to
a point $B\in\cA_{g-1}$ in the Satake compactification, and work in
the universal level family over it. As before, we denote by $n=C\cdot
\gT_B$; the branch divisor on a test curve is just a number of
points, which we thus compute in terms of $n$.

Indeed, by the first part of the proof of \ref{Pmumbound} we have,
substituting the expression for the class from Theorem \ref{Ttheta}:
$$
 \BR(f)=\sum_{i=0}^{m-1}\calD_i^{g+1}=\sum_{i=0}^{m-1}(\xi_i+m\ti{\mu_i}
 +\frac12\ti{\ga'}_i+\frac{m}4\ti{\eta}_i)^{g+1}.
$$
If we now denote $\mu_i':=m\mu_i$ and $\eta_i':=m\eta_i$ (and recall
that $\ga_i'=m\ga_i$), then all of the identities that were valid for
$\xi,\mu,\eta,\ga$ in the no level case are valid for
$\xi_i,\mu_i',\eta_i',\ga_i$ in the level case, except that now we
get an extra factor of $m^g$ for the top intersection numbers, i.e.
we have
\[
  \begin{aligned}
    (\blacksquare')&\quad \xi_i^2=-\ti{\ga_i'}\xi_i\\
  (\square')& \quad  (\eta'_i)^2=0
  \qquad\qquad&(\diamondsuit')\quad& (\eta_i')(\mu_i')^{g-1}&=&\ m^gn(g-1)!\\
  (\triangle')&\quad \ga_i'\eta_i'=0
  \qquad\qquad&(\heartsuit')\quad& (\ga_i')^k(\mu_i')^{g-k}&=&
        \ \begin{cases}-2m^g(g-2)!{\rm\ if\ } k=2\\
                      0{\rm\ otherwise}
        \end{cases}
  \end{aligned}
\]
It thus follows that if we use these classes in identity $(T)$, there
is simply an extra factor of $m^g$, and thus the computation in
\ref{Pmumbound} carries over verbatim to yield
$$
\BR(f)=\sum_{i=0}^{m-1}\calD_i^{g+1}=m\left(m^g\frac{n(g+1)!}{6}\right)
$$
as expected.
\end{exm}


\end{document}